\documentclass[12pt]{article}
\usepackage{amsmath,amssymb,amsthm,amscd,a4wide}
\def\be{\begin{equation}}\def\ee{\end{equation}}
\def\lb{\label}
\begin{document}
\begin{center} {\Large \bf Idempotents for Birman--Murakami--Wenzl algebras \\ \vspace{.3cm}
and reflection equation}

\vspace{1.5cm}
{\large  A. P. Isaev$^{\diamond}$, A. I. Molev$^{\circ}$ and O. V. Ogievetsky$^{\ast}$\footnote{On leave of absence
from P. N. Lebedev Physical Institute, Leninsky Pr. 53,
117924 Moscow, Russia}}

\vspace{1cm}
$\diamond\ ${Bogoliubov Laboratory of Theoretical Physics\\ Joint Institute for Nuclear Research\\
Dubna, Moscow region 141980, Russia\\ isaevap@theor.jinr.ru}

\vspace{.6cm}
$\circ\ ${School of Mathematics and Statistics\\ University of Sydney,
NSW 2006, Australia\\ alexander.molev@sydney.edu.au}

\vspace{.6cm}
$\ast\ ${Center of Theoretical Physics\footnote{Unit\'e Mixte de Recherche (UMR 6207) du CNRS et des
Universit\'es Aix--Marseille I, Aix--Marseille II et du Sud Toulon -- Var; laboratoire
affili\'e \`a la FRUMAM (FR 2291)}, Luminy \\ 13288 Marseille, France\\
\vspace{.2cm}
and\\
\vspace{.2cm}
J.-V. Poncelet French-Russian Laboratory, UMI 2615 du CNRS\\
Independent University of Moscow, 11 B. Vlasievski per., 119002 Moscow, Russia\\
oleg@cpt.univ-mrs.fr}
\end{center}

\vskip 1cm
\begin{abstract}

\vskip .3cm
\noindent
A complete system of  pairwise orthogonal
minimal idempotents for Birman--Mura\-ka\-mi--Wenzl algebras is obtained by a consecutive evaluation of
a rational function in several variables  on sequences of quantum
contents of up-down tableaux. A by-product of the construction is a one-parameter family of fusion procedures for
Hecke
algebras. Classical limits to two different fusion procedures for Brauer algebras are described.
\end{abstract}

\vspace{3cm}

\section{Introduction}

Let $V$ be a vector space equipped with a non-degenerate bilinear (symmetric or antisymmetric) form.
Let $G$ be the group of linear transformations of $V$ preserving the form.
The Brauer--Schur--Weyl duality relates tensor representations of $G$ with representations of the Brauer algebras.
The ``quantum group" version of the Brauer--Schur--Weyl duality relates tensor representations of a
quantum group associated with $G$ with representations of the Birman--Murakami--Wenzl (BMW) algebras, introduced in
\cite{BirWen} and \cite{Mur0}. A quantum group associated with $G$ can be defined through an $R$-matrix of the
``$BCD$-type";
it is an endomorphism of $V\otimes V$, having three eigenvalues with prescribed multiplicities, which
satisfies the Yang--Baxter equation together with extra conditions ensuring that the tensor powers of $V$
carry a ``local" representation of the tower of the BMW algebras. Thus,
in contrast to the classical situation,
the quantum version of the Brauer--Schur--Weyl
duality requires a choice of a $BCD$-type $R$-matrix.
It should be noted that
the problem  of a description of $BCD$-type $R$-matrices
is difficult and far from being completed (it is an open question already for $R$-matrices of the $A$-type, see
\cite{O}). To
use the quantum version of the Brauer--Schur--Weyl duality one needs to
know explicit bases in the representation spaces of the BMW algebras. This task splits into two parts. First, one
needs
to know a complete set of pairwise orthogonal minimal idempotents of the BMW algebras; second, to calculate, given
a $BCD$-type $R$-matrix, the image of these idempotents in the corresponding local representation. The present paper
is a contribution to the first, $R$-matrix independent, part of the task, it deals with formulas for minimal
idempotents of the BMW
algebras. Such formulas can be written down in terms of the generators of the BMW algebras,
without a concrete choice of an $R$-matrix. This requires a good understanding of the representation theory of the
BMW
algebras. Different aspects of the representation theory of the BMW algebras were discussed in
the literature, see {\it e.g.} \cite{BeBla}, \cite{BirWen}, \cite{LeRa}, \cite{RaWe} and \cite{We}. The family of the
BMW-algebras
depends on a discrete parameter $n=0,1,2,\dots$ and two continuos parameters, $q$ and $\nu$.
 Formulas for the action of the generators of the
BMW-algebras on spaces of their irreducible representations were given in  \cite{Mur0}.
 The eigenvalues of the so called Jucys--Murphy elements were calculated
in \cite{LeRa}, see also \cite{IsOg} (the article \cite{LeRa} does not use terminology ``Jucys--Murphy elements"). One
can show that the Jucys--Murphy elements generate the maximal
commutative subalgebra of the BMW algebra; in particular, the Jucys--Murphy elements distinguish the basis vectors and
thus imply formulas for
the complete set of pairwise orthogonal minimal idempotents.

\vskip .2cm

We suggest another formula for the same set of idempotents in the present paper.
Our way to describe the idempotents is in the spirit of what is often referred to as the {\it fusion procedure}. For
generic
values of the parameters of the BMW algebra, basis vectors in the irreducible representations are in one-to-one
correspondence with sequences of {\it quantum contents}, see precise definitions in Section \ref{secId}. We introduce
in Section  \ref{secId} a certain BMW algebra-valued rational function in several variables (which we call {\it
fusion} function).
The idempotent, corresponding to a given content sequence, is then obtained by a
consecutive evaluation of the fusion function at values of these variables   given by the quantum contents.
The fusion procedure (for the symmetric group) originates in the work of Jucys \cite{Ju0}, see also \cite{Cher1} and
\cite{Naz}. A version of the fusion procedure for the symmetric group involving the consecutive evaluation
appeared in \cite{Mol}. The fusion function was written down for the Hecke algebras and Brauer algebras, see
\cite{IsMol}, \cite {IMO}, \cite{IsMoOs} and \cite{Naz2}.

\vskip .2cm
The BMW algebra admits a natural quotient isomorphic to the Hecke algebra.
The BMW algebra depends on a certain parameter $\nu$, but this dependence
disappears in the Hecke algebra quotient. Nevertheless our fusion function for the BMW algebra
depends on $\nu$ non-trivially. As a consequence, the image of the BMW fusion function in the Hecke algebra
quotient provides a one parameter family, see paragraph 5 of Section \ref{secId}, of the fusion functions for the
Hecke
algebras.

\vskip .2cm
In paragraph 6 Section \ref{secId} we make an observation on a deep relation between the BMW fusion function
and the solution, presented in \cite{IsOg2}, of the so called reflection equation. Thus for algebras related to the
centralizer construction for classical and quantum groups of types $GL$, $O$ and $Sp$ (and $OSp$ for super-groups)
in the defining representation (like the group algebra of the symmetric group, the Hecke algebra, the Brauer algebra
and the Birman--Murakami--Wenzl algebra) the fusion function is built with the help of solutions of the reflection
equation.

\vskip .2cm
The classical version of the reflection equation was used in \cite{IMO} to construct certain evaluation
homomorphisms.
We postpone the study of the evaluation homomorphisms for the quantum algebras for a future work.

\vskip .2cm
Two different versions of the fusion functions for the Brauer algebra are given in \cite{IsMol} and \cite {IMO}. The
Brauer
algebra can be obtained in several ways as a contraction (or classical limit) of the BMW algebra. In Section
\ref{secclli}
we derive both fusion functions for the Brauer algebra from the same fusion function for the BMW algebra using two
different
contractions.

\section{Birman--Murakami--Wenzl algebra}
\setcounter{equation}{0}

\paragraph{1. Definition and basic relations.} \hspace{-.2cm}The {\it Birman--Murakami--Wenzl algebra}
$BMW_n(q,\nu)$ was defined in \cite{BirWen}, \cite{Mur00}
and \cite{Mur0}. It is generated over $\mathbb{C}$ by invertible
elements $T_1,\dots,T_{n-1}$ with the following defining relations
\be\lb{bmw01}T_i\, T_{i+1}\, T_i=T_{i+1}\, T_i\, T_{i+1}\; ,\;\;\; T_i\, T_{j}=T_j\, T_{i}\;\;\; {\rm for}\;\; |i-j|
>1\; ,\ee
\be\lb{bmw1}\kappa_i\, T_i=T_i\,\kappa_i=\nu\,\kappa_i\; ,\phantom{T_j}\ee
\be\lb{bmw2}\kappa_i\, T_{i-1}^{\varepsilon}\,\kappa_i=\nu^{-\varepsilon}\,\kappa_i\;\; ,\;\;
\kappa_i\, T_{i+1}^{\varepsilon}\,\kappa_i=\nu^{-\varepsilon}\,\kappa_i\ \ \text{with}\ \varepsilon=\pm 1\; ,\ee
where
\be\lb{bmw3}\kappa_i:=1-\frac{T_i-T_i^{-1}}{q -q^{-1}}\; .\ee
Here $q$ and $\nu$ are complex parameters of the algebra which we assume generic in the sequel; in particular, the
definition
(\ref{bmw3}) makes sense, the denominator in the right hand side does not vanish.

\vskip .2cm
The quotient of the algebra $BMW_n(q,\nu)$ by the ideal generated by the elements $\kappa_1,\dots,\kappa_n$
(in fact, this ideal is generated by any one of these elements, say, $\kappa_1$) is isomorphic to the Hecke algebra
$H_n(q)$. We shall often omit the parameters in the notation for the algebras and write simply $BMW_n$
and $H_n$.

\vskip .2cm
Let
\be\label{lammu}\mu =\displaystyle{\frac{q -q^{-1}+ \nu^{-1} - \nu}{q -q^{-1}}}=
\displaystyle{\frac{(q^{-1}+\nu)(q-\nu)}{\nu(q -q^{-1})}}\ .\ee
The following relations can be derived from  (\ref{bmw01})--(\ref{bmw2}):
\begin{equation}\lb{bmw4}\kappa_i^2 =  \mu \, \kappa_i
\end{equation}
then, with $\varepsilon =\pm 1$,
\begin{eqnarray}
\lb{bmw5}\kappa_i\, T_{i+\varepsilon}\, T_i\!\!  &=&\!\!  T_{i+\varepsilon}\, T_i \,\kappa_{i+\varepsilon} \;
,\\[.5em]
\lb{bmw8}\kappa_i\,\kappa_{i+\varepsilon}\,\kappa_i\!\!  &=&\! \! \kappa_{i}\; ,\\[.5em]
\lb{bmw9}\bigl(T_i-\!(q -q^{-1})\bigr)\kappa_{i+\varepsilon}\bigl(T_i-\!(q -q^{-1})\bigr)\!\!   &=&\!\!
\bigl( T_{i+\varepsilon}-\!(q -q^{-1})\bigr)\kappa_{i}\bigl(T_{i+\varepsilon} -\!(q -q^{-1})\bigr)\, ,\\[.5em]
\lb{bmw8a}T_{i+\varepsilon}\,\kappa_i\, T_{i +\varepsilon}\!\!  &=&\!\!  T^{-1}_i\,\kappa_{i+\varepsilon}\, T^{-1}_i\;
,
\end{eqnarray}
and
\begin{eqnarray}
\lb{bmw6}\kappa_i\, T_{i+\varepsilon}\, T_{i}&=&\kappa_i\,\kappa_{i+\varepsilon}\; ,\\[.5em]
\lb{bmw7}\kappa_i\, T^{-1}_{i+\varepsilon}\, T^{-1}_{i}&=&\kappa_i\,\kappa_{i+\varepsilon}\; ,\\[.5em]
\lb{bmw8b}\kappa_{i+\varepsilon}\,\kappa_i\,\bigl(T_{i+\varepsilon}-(q -q^{-1})\bigr) &=&
\kappa_{i+\varepsilon}\,\bigl( T_i-(q -q^{-1})\bigr)\; ,\end{eqnarray}
together with their images under the anti-automorphism $\rho$ of the algebra $BMW_n$ defined on the generators by
\be\lb{aau} \rho(T_i)=T_i\; .\ee

\paragraph{2. Baxterized elements.} The {\it baxterized elements} $T_i(u,v)$ are defined by
\be\lb{a00}T_i(u,v):=T_i+\frac{q-q^{-1}}{v/u-1}+\frac{q-q^{-1}}{1+\nu^{-1}qv/ u} \kappa_i\; ,\ee
see  \cite{CGX}, \cite{Isa}, \cite{Jones} and \cite{Mur}. They are rational functions in complex variables $u$ and $v$
which
are called {\it spectral variables}. The elements $T_i(u,v)$ depend on the ratio of the spectral parameters;
later we shall evaluate spectral variables $u$ and $v$ in terms of contents of boxes of Young diagram; this explains
the use
of having both spectral variables in the notation (\ref{a00}) for the baxterized element. However for brevity we
shall
denote sometimes the baxterized elements by $T_i(u/v)$ (with one argument only), $T_i(u):=T_i(u,1)$.

\vskip .2cm
The baxterized elements satisfy the braid relation of the form
\be\lb{ybe0}T_i(u_{2},u_{3})T_{i+1}(u_1,u_{3})T_i(u_1,u_{2})=T_{i+1}(u_1,u_{2})T_i(u_1,u_{3})T_{i+1}(u_{2},u_{3})\;
.\ee
The inverses of the baxterized elements are given by
\be\lb{aa}T_i(v,u)^{-1}=T_i(u,v)\, f(u,v)\; ,\ee
where
\be f(u,v)=\frac{(u-v)^2}{(u-q^2 v)(u-q^{-2}v)}=f(v,u)\ .\ee
Let
\be\lb{a0}Q_i(u,v;c):=T_i \left(\frac{1}{cuv}
\right)=T_i+\frac{q-q^{-1}}{cuv-1}+\frac{q-q^{-1}}{1+\nu^{-1}qcuv}\kappa_i \; .\ee
We keep three arguments $c,u$ and $v$ for the later convenience.

\vskip .2cm
It follows from (\ref{ybe0}) that
\be\lb{tqq}T_i(u_2,u_3)Q_{i+1}(u_1,u_3;c)Q_i(u_1,u_2;c)=Q_{i+1}(u_1,u_2;c)Q_i(u_1,u_3;c)T_{i+1}(u_2,u_3)\; .\ee

\paragraph{3. Jucys--Murphy elements.}
The Jucys--Murphy elements of the algebra $B\!M\!W_n$ are defined by
\be\lb{jme}y_1=1 \;, \;\;\; y_{k+1}=T_k\dots T_2\, T_1^2\, T_2\dots T_k\ , \ \ k=1, \dots ,n-1\; .\ee
The elements $y_1,\dots,y_n$ pairwise commute and satisfy the identities
\be\lb{yyk}\kappa_j\, y_{j+1}\, y_j=y_j\, y_{j+1}\,\kappa_j=\nu^2\,\kappa_j\; .\ee

The Jucys–Murphy elements were originally used for constructing idempotents for the symmetric groups in [18], [27].
Analogues of the Jucys–Murphy elements can be defined for a number of important algebras related to the symmetric
group rings; they turn out to generate maximal commutative subalgebras in these rings (see \cite{IMO}, \cite{IsOg3}, \cite{OPdA},
\cite{OV} and references therein).
The commutative subalgebra, generated by the Jucys--Murphy elements $y_1,\dots,y_n$, of the generic
algebra $B\!M\!W_n$  is maximal as well; it follows from the results in \cite{IsOg},\cite{LeRa}.
In the next section we use this maximality for the construction of the idempotents
for the BMW algebra.

\section{Idempotents for BMW algebra}\label{secId}
\setcounter{equation}{0}

\paragraph{1. Up-down tableaux.} Recall that the Bratteli diagram for the chain of generic BMW algebras is the
up-down Young graph;
thus the irreducible representations of the generic algebra $BMW_{n}$ are in one-to-one correspondence with the
partitions
$\lambda$ of $n-2k$, $0\leqslant k\leqslant \lfloor n/2\rfloor$. We use the convention
to depict a partition $\lambda =(\lambda_1,\lambda_2,\dots)$ by its Young diagram which is the left-justified array of
rows of
boxes containing $\lambda_1$ boxes in the first row,  $\lambda_2$ boxes in the second row and so on. Basis vectors
$v_{_{U_n}}$ in the irreducible representation
corresponding to $\lambda$ are indexed by the up-down tableaux of shape $\lambda$, that is, by sequences
$U_n=(\Lambda_1,\Lambda_2,\dots,\Lambda_n)$ of Young diagrams, such that $\Lambda_1=\Box$, $\Lambda_j$ differs
from $\Lambda_{j-1}$ by exactly one box (added or removed), $j=2,\dots,n$, and $\Lambda_n=\lambda$; see \cite{GHJ}
for necessary definitions of Bratteli diagrams {\it etc}.

\vskip .2cm
The usual notion of transposition of Young diagrams and tableaux extends to up-down tableaux;  for an up-down
tableau $U_n=(\Lambda_1,\Lambda_2,\dots,\Lambda_n)$ the sequence
$(\Lambda_1^t,\Lambda_2^t,\dots,\Lambda_n^t)$, formed by the transposed Young diagrams,
is clearly an up-down tableau and we set $U_n^t:=(\Lambda_1^t,\Lambda_2^t,\dots,\Lambda_n^t)$.

\vskip .2cm
The up-down Young graph arises in the study of the common spectrum of the Jucys--Murphy elements within the algebra
itself
\cite{IsOg}; as for representations, the images of the Jucys--Murphy elements are diagonal in the chosen basis,
\be\label{jmevobv}\rho_{\lambda}(y_j)v_{_{U_n}}=c_j(U_n)v_{_{U_n}}\ .\ee
Here $v_{_{U_n}}$ is the basis vector associated to the up-down tableau
$U_n=(\Lambda_1,\Lambda_2,\dots,\Lambda_n)$ of shape
$\lambda$ and $\rho_{\lambda}(z)$ is the operator of the element $z\in BMW_n$ in the representation corresponding to
$\lambda$;  the number $c_j(U_n)$ is the quantum content of the box added/removed at the step $j$ of building
the up-down tableau $U_n$; for the box in the $a$-th row and $b$-th column the quantum content is defined by
\be\label{doqcfudt}c_j(U_n)=q^{2(b-a)}\ee
if this box was added to the diagram $\Lambda_{j-1}$ to obtain the diagram $\Lambda_j$ and
\be\label{doqcfudtb} c_j(U_n)=\nu^2 q^{2(-b+a)}\ \ee
if this box was removed from the diagram $\Lambda_{j-1}$ to obtain the diagram $\Lambda_j$.
The spectrum of the Jucys--Murphy elements in the irreducible representations of the BMW algebra
was calculated in \cite{LeRa}.

\vskip .2cm

Note that
$c_1=1$. In the generic regime, the up-down tableau $U_n$ is uniquely determined by its content sequence
$\{ c_1(U_n),\dots,c_n(U_n)\}$.

\vskip .2cm
We denote by $E_{U_n}\equiv E^{\lambda}_{U_n}\in BMW_n$ the primitive idempotent related to the vector $v_{U_n}$ in
the
irreducible representation corresponding to a diagram $\lambda$. That is, $E_{U_n}$ is the projector onto the
1-dimensional subspace
spanned by the vector $v_{U_n}$ along the subspace spanned by the other basis vectors, and $E_{U_n}$  annihilates
all other irreducible representations.

\paragraph{2. Inductive formula for idempotents of BMW algebras.} {}For an up-down tableau
$U_n=(\Lambda_1,\Lambda_2,\dots,\Lambda_{n-1},\Lambda_n)$ of
length $n$ we denote by $U_k$, $k\leqslant n$, the initial segment of $U_n$ of length $k$,
$U_k:=(\Lambda_1,\Lambda_2,\dots,
\Lambda_k)$. {}For brevity, let $c_j=c_j(U_n)$ and let $c$ be a parameter, $c\neq c_n^{-2}$.
Note that by (\ref{doqcfudt}) the Jucys--Murphy element $y_n$ takes the value $c_n$ on the image of the idempotent
$E_{U_n}$. Suppose that $Y$ is an eigenvalue of $y_n$ on the image of $E_{U_{n-1}}$. {}For $Y\neq c_n$ the rational function
$\displaystyle{\frac{(cu\, Y-1)}{(u-Y)}\frac{(u-c_n)}{(cuc_n-1)}}$ vanishes at  $u=c_n$ while
for $Y=c_n$ this rational function
simplifies to 1. Thus,
\be\lb{bran}E_{U_n}=\left. E_{U_{n-1}} \frac{(cu\, y_n-1)}{(u-y_n)}\frac{(u-c_n)}{(cuc_n-1)}\right|_{u=c_n}\; \ee
with arbitrary $c\neq c_n^{-2}$. The parameter $c$ will be fixed below:
\be\lb{ccc}c=-q^{-1}\nu^{-1} \; .\ee
The form of the rational function in the right hand side of (\ref{bran}) and
the choice (\ref{ccc}) of $c$ are suggested by the proof of Lemma 1. With the choice (\ref{ccc}) the inequality $c\neq
c_n^{-2}$ holds in the
generic regime.

\paragraph{3. Preparatory lemma.} We shall rewrite the product $\displaystyle{E_{U_{n-1}}\hspace{.05cm}\frac{cu\,
y_n-1}{u-y_n}}$ in the right hand side
of (\ref{bran}). Define a family of $BMW_n$-valued rational functions by
\be\label{a40b}{\cal Y}_1(u):=\displaystyle{\frac{cu\, y_1-1}{u-y_1}}\ee
and, for $j=2,\dots,n$,
\be\lb{a4}\hspace{-.4cm}\begin{array}{c}{\cal Y}_j(u_1,\dots ,u_j):=Q_{j-1}(u_{j-1},u_j;c)\,
{\cal Y}_{j-1}(u_1,\dots ,u_{j-2},u_j)\, T_{j-1}(u_j,u_{j-1})^{-1}\\ [1em]
\hspace{1.5cm}=Q_{j-1}(u_{j-1},u_j;c)\dots Q_1(u_1,u_j;c) \; {\cal Y}_1(u_j)\;T_1^{-1}(u_j,u_1)\dots
T_{j-1}(u_j,u_{j-1})^{-1} \; .\end{array} \ee

\vskip .3cm
\noindent {\bf Lemma 1.}
\be\lb{a5}E_{U_{n-1}}\, {\cal Y}_n(c_1,c_2,\dots ,c_{n-1},u)=E_{U_{n-1}} \frac{cu \, y_n -1}{u- y_n} \; .\ee

\vskip .2cm\noindent
{\it Proof.}  We prove (\ref{a5}) by induction on $n$. The induction base is $n=1$ ($E_{U_0}=1$, $U_0$ is
empty) is tautological.
By (\ref{a4}) and the induction hypothesis, (\ref{a5}) reduces to the equality
\be\lb{a6}E_{U_{n-1}}\, Q_{n-1}(c_{n-1},u;c)\,\frac{cu\, y_{n-1} -1}{u- y_{n-1}}\,T_{n-1}(u,c_{n-1})^{-1}=
E_{U_{n-1}} \frac{cu \, y_n -1}{u- y_n} \; .\ee
Since $y_n$ commutes with $E_{U_{n-1}}$, this can be rewritten in the form
\be\lb{a8}E_{U_{n-1}} (u- y_n)Q_{n-1}(c_{n-1},u;c)(cu\, y_{n-1}-1)\!
=\!E_{U_{n-1}} (cu \, y_n-1)T_{n-1} (u,c_{n-1}) (u- y_{n-1}),\ee
or, with the definition (\ref{a0}),
\be\lb{a9}\hspace{-.05cm}\begin{array}{c}E_{U_{n-1}}\,  (u- y_n) \bigl(T_{n-1}+\displaystyle{\frac{q-q^{-1}}{c
c_{n-1}u -1}+
\frac{q-q^{-1}}{1+\nu^{-1}qcc_{n-1}u}} \kappa_{n-1}\bigr) (cu \, y_{n-1} -1)  \\[1.5em]
\hspace{1.5cm}= E_{U_{n-1}} (cu\, y_n -1)\bigl(T_{n-1}+\displaystyle{\frac{(q-q^{-1})u}{c_{n-1}
-u}+\frac{(q-q^{-1})u}{u+
\nu^{-1}qc_{n-1}}} \kappa_{n-1}\bigr)(u- y_{n-1}) \; .\end{array}\ee
We transform the left hand side of (\ref{a9}) to
\be\lb{a10}\begin{array}{c}E_{U_{n-1}}\, (u-y_n)\Bigl(cu y_n\bigl(\underline{T_{n-1}-(q-q^{-1})}+(q-q^{-1})
\kappa_{n-1}\bigr)-\underline{T_{n-1}}  \\[1.5em]
+\displaystyle{\frac{q-q^{-1}}{c c_{n-1}u -1}\underline{(cu\,
c_{n-1}-1)}+\frac{q-q^{-1}}{1+\nu^{-1}qcc_{n-1}u}}\kappa_{n-1}
(cu \, y_{n-1} -1)\Bigr)\end{array}\ee
and the right hand side of (\ref{a9}) to
\be\lb{a11}\begin{array}{c}E_{U_{n-1}}\, (cuy_n-1)\Bigl(-y_n\bigl(\underline{T_{n-1}-(q-q^{-1})}+(q-q^{-1})
\kappa_{n-1}\bigr)+\underline{uT_{n-1}}  \\[1.5em]
+\displaystyle{\frac{(q-q^{-1})u}{c_{n-1} -u} \underline{(u- c_{n-1})}+\frac{(q-q^{-1})u}{u+\nu^{-1} q c_{n-1}}}
\kappa_{n-1}
(u-y_{n-1})\Bigr)\ .\end{array}\ee
The underlined terms in (\ref{a10}) cancel the underlined terms in (\ref{a11}) and we are left to verify that
\be\lb{a12}\begin{array}{c}E_{U_{n-1}} \,  (u- y_n) \Bigl(cu
y_n\kappa_{n-1}+\displaystyle{\frac{1}{1+\nu^{-1}qcc_{n-1}u}}
\kappa_{n-1}(cu \, y_{n-1} -1)\biggr) \\[1em]
\hspace{1.5cm}=E_{U_{n-1}} \,  (cuy_n-1)
\biggl(-y_n\kappa_{n-1}+\displaystyle{\frac{u}{u+\nu^{-1}qc_{n-1}}}\kappa_{n-1}
(u- y_{n-1})\Bigr) \; ,\end{array}\ee
where, recall, $c=-q^{-1}\nu^{1}$.
We use the relations $E_{U_{n-1}} y_{n-1} = E_{U_{n-1}} c_{n-1}$ and (\ref{yyk}) which imply
\be\lb{a13}E_{U_{n-1}}
y_n\kappa_{n-1}=E_{U_{n-1}} \kappa_{n-1}\frac{\nu^2}{c_{n-1}}\ .\ee
Now the check of (\ref{a12}) is straightforward. \hfill$\Box$

\paragraph{4. Fusion procedure.} Define
the {\it fusion} function $\Psi_{U_n}$ by
\be\lb{abi00}\Psi_{U_n}(u_1,u_2,\dots ,u_n):= \left( \prod\limits_{k=1}^n \frac{u_k-c_k}{c u_k c_k -1} \right)
{\cal Y}_1(u_1) {\cal Y}_2(u_1,u_2) \dots {\cal Y}_n(u_1, \dots ,u_n) \; .\ee
Sometimes we shall write, instead of the label $U_n$ of the fusion function, the corresponding sequence of contents.

\vskip .2cm
The following theorem is the main result of the paper.

\vspace{0.2cm}
\noindent {\bf Theorem 2.} {\it The idempotent $E_{U_n}$ is equal to the  subsequent evaluation of
$\Psi_{U_n}$ at the points $u_i = c_i$,
\be\lb{idemp}E_{U_n}=\left.\left. \Psi_{U_n}(u_1,u_2,\dots ,u_n)  \right|_{u_1=c_1} \cdots\right|_{u_n=c_n}\; .\ee}

\vspace{0.1cm}
\noindent {\it Proof.} We prove the equality (\ref{idemp}) again by induction on $n$. The base of induction, $n=1$,
is
immediate (recall that $c_1=1$):
$$\Psi_{\Box}(u_1)|_{u_1=c_1}= \frac{u_1-c_1}{c u_1 c_1 -1} {\cal Y}_1(u_1) |_{u_1=c_1}= \left.
\frac{u_1-c_1}{c u_1 c_1 -1}\, \frac{cu_1\, c_1-1}{u_1-c_1} \right|_{u_1=c_1} =1=E_{U_1}  \; .$$
Assume that (\ref{idemp}) is valid for the tableaux $U_{n-1}$. Then
\be\lb{a55}\begin{array}{c}\left. \left.  \Psi_{U_n}(u_1,u_2,\dots ,u_{n-1},u)\right|_{u_{1}=c_1} \cdots
\right|_{u_{n-1}=c_{n-1}}=
\displaystyle{\frac{u-c_n}{cu\, c_n -1}} E_{U_{n-1}} {\cal Y}_n(c_1, \dots ,c_{n-1},u)\\[1.5em]
=\displaystyle{\frac{u- c_n}{cu \, c_n -1} E_{U_{n-1}} \frac{cu \, y_n -1}{u- y_n}}\; ;\end{array}\ee
we used the induction hypothesis in the first equality and Lemma 1 in the second equality.
By (\ref{bran}), the evaluation of the last expression in (\ref{a55}) at $u=c_n$ is equal to $E_{U_n}$.\hfill $\Box$

\paragraph{Example 1.} We shall illustrate Theorem 2 on the simplest example of the algebra $BMW_2$.
The algebra $BMW_2$ is commutative; it is generated by $T_1$ which has the projector decomposition
\be\label{prdotg}T_1=qS-q^{-1}A+\nu\Pi\ . \ee
The projector
$S$ is called symmetrizer - its image in the Hecke algebra quotient tends to the symmetrizer in the classical limit
$q\rightarrow 1$; similarly, the projector $A$ is called antisymmetrizer. One has
\be\label{idintots}S=\frac{(T+q^{-1})(T-\nu)}{(q+q^{-1})(q-\nu)}=\frac{1}{q+q^{-1}}\Bigl(T_1+q^{-1}+
\frac{q-q^{-1}}{1- q \nu^{-1}} \kappa_1\Bigr)\ ,\ee
\vskip .1cm
\be\label{idintota} A=\frac{(T-q)(T-\nu)}{(-q^{-1}-q)(-q^{-1}-\nu)}=-\frac{1}{q+q^{-1}}\Bigl(T_1-q+
\frac{q-q^{-1}}{1+q^{-1}\nu^{-1}} \kappa_1\Bigr)\ ,\ee
\vskip .2cm
\be\label{idintotp} \Pi=\frac{(T-q)(T+q^{-1})}{(\nu-q)(\nu +q^{-1})}= \frac{1}{\mu} \kappa_1\; ,\ee
where $\mu$ is given by (\ref{lammu}). The elements $S$, $A$ and $\Pi$ form a complete system of pairwise orthogonal
primitive idempotents of $BMW_2$. Every element $x$ of $BMW_2$ is a linear combination of  $S$, $A$ and $\Pi$; the
coefficients can be found by calculating the products of $x$ with these three idempotents.

\vskip .2cm
The algebra $BMW_2$ has three irreducible representations, all three are one-dimensional.
We have three up-down tableaux with content sequences $\{ 1,q^2\}$, $\{ 1,q^{-2}\}$ and $\{ 1,\nu^2\}$.
The fusion function $\Psi_{\{ c_1,c_2\} }$ reads
$$\Psi_{\{ c_1,c_2\} }(u_1,u_2) = \left( \prod\limits_{k=1}^2 \frac{u_k-c_k}{c u_k c_k -1} \right)
\frac{cu_1\, y_1 -1}{u_1-y_1}  \,  Q_1(u_1,u_2;c) \frac{cu_2\, y_1 -1}{u_2-y_1} T_1(u_2,u_1)^{-1}$$
\be\lb{a1}=
\frac{u_2-c_2}{c u_2 c_2 -1}\,\frac{cu_2-1}{u_2-1}
\frac{(u_1 - u_2)^2}{(u_2 - q^2 u_1)(u_2 - q^{-2} u_1)}T_1 \Bigl(\frac{1}{c u_1 u_2}\Bigr) \,
T_1\Bigl(\frac{u_1}{u_2}\Bigr) \; .\ee
We used the fact that $y_1=c_1=1$. Below we write $E_{\{ c_1,c_2\} }$ for the idempotent $E_{U_2}$ corresponding to
the
up-down tableau $U_2$ with the content sequence $\{ c_1,c_2\}$.

\vskip .2cm
\noindent 1) Content sequence $\{1,q^2\}$. We have
\be\label{udt2s}E_{\{ 1,q^2\} }=\Psi_{\{ 1,q^2\} }(u_1,u_2)
|_{u_1=1}|_{u_2=q^2}=\frac{q\nu^{-1}+1}{q^3\nu^{-1}+1}\frac{q}{q+q^{-1}}
T_1(-q^{-1}\nu)T_1(q^{-2})\ .\ee
One finds $T_1(-q^{-1}\nu)\Pi =0$ and $T_1(q^{-2})A=0$ so the last expression in (\ref{udt2s}) is proportional to the
symmetrizer. Calculating $T_1(-q^{-1}\nu)S =q^{-1}\displaystyle{\frac{q^3\nu^{-1}+1}{q\nu^{-1}+1}}S$ and
$T_1(q^{-2})S=(q+q^{-1})S$ we obtain
$$E_{\{ 1,q^2\} }=S\ .$$

\noindent 2) Content sequence $\{ 1,q^{-2}\}$. We have
\be\label{udt2a}E_{\{ 1,q^{-2}\} }=\Psi_{\{ 1,q^{-2}\} }(u_1,u_2)
|_{u_1=1}|_{u_2=q^{-2}}=\frac{q^{-3}\nu^{-1}+1}{q^{-5}\nu^{-1}+1}
\frac{q^{-1}}{q+q^{-1}}
T_1(-q^3\nu)T_1(q^{2})\ .\ee
This time the same factor $T_1(q^{2})$ is annihilated by two idempotents, $T_1(q^{2})S=0$ and $T_1(q^{2})\Pi =0$;
in fact, $T_1(q^2)= -(q+q^{-1})A$. Calculating $T_1(-q^3\nu)A
=-q\displaystyle{\frac{q^{-5}\nu^{-1}+1}{q^{-3}\nu^{-1}+1}}A$
we obtain
$$E_{\{ 1,q^{-2}\} }=A=-\frac{1}{q+q^{-1}}T_1(q^2)\ .$$

\vskip .2cm
\noindent 3) Content sequence $\{ 1,\nu^2\}$. We have
\be\label{udt2p}\begin{array}{rcl}E_{\{ 1,\nu^2\} }&=&\Psi_{\{ 1,\nu^2\} }(u_1,u_2) |_{u_1=1}|_{u_2=\nu^2}\\[1em]
&=&\displaystyle{\frac{(u_2-\nu^2)(u_2-1)(q^{-1}\nu^{-1}u_2+1)}{(u_2-q^2)(u_2-q^{-2})(q^{-1}\nu u_2+1)}}
T_1(-q\nu u_2^{-1})T_1(u_2^{-1})|_{u_2=\nu^2}\ .\end{array}\ee
Here the scalar coefficient has a zero at $u_2=\nu^2$, this is a new phenomenon compared to (\ref{udt2s})
and (\ref{udt2a}). The factor $T_1(u_2^{-1})$ is regular at $u_2=\nu^2$.
In the factor $T_1(-q\nu u_2^{-1})$ the coefficient in front of $\kappa_1$ is
$\displaystyle{\frac{q-q^{-1}}{1-\nu^{-2}u_2}}$
and has a simple pole at $u_2=\nu^2$; other coefficients are regular at $u_2=\nu^2$. Thus, only $\Pi$ survives in
 $T_1(-q\nu u_2^{-1})$. Calculating $T_1(\nu^2)\Pi =\displaystyle{\frac{\nu^{-1}q(q^{-1}\nu^3+1)(\nu
 -q^{-1})}{\nu^2-1}}\Pi$
we obtain
$$E_{\{ 1,\nu^2\} }=\Pi\ .$$

\paragraph{Example 2.} One can rewrite the expression (\ref{abi00}) for the fusion function in the form
\be\lb{abi00b}\Psi_{U_n}(u_1,u_2,\dots ,u_n)=\mathfrak{y}_{U_n}(u_1,u_2,\dots ,u_n)Y_n(u_1,u_2,\dots ,u_n)\ ,\ee
where
\be\label{abi00c}\mathfrak{y}_{U_n}(u_1,u_2,\dots ,u_n):=\prod\limits_{k=1}^n\left(\frac{u_k-c_k}{cu_kc_k-1}\;
\frac{cu_k-1}{u_k-1}
\prod\limits_{r=1}^{k-1} \frac{(u_k-u_r)^2}{(u_k-q^2u_r)(u_k-q^{-2}u_r)}\right)\ee
and
\be\label{abi00d}Y_n(u_1,u_2,\dots ,u_n):=
{\cal Q}_2(u_1,u_2)\dots {\cal Q}_n(u_1, \dots ,u_n)\ {\cal T}_n(u_1,\dots,u_n)
\dots {\cal T}_2(u_1,u_2)  \; \ee
with
\be\label{abi00e}{\cal Q}_j(u_1, \dots ,u_j):=T_{j-1}\Bigl(\frac{1}{cu_1u_j}\Bigr)\dots
T_1\Bigl(\frac{1}{cu_{j-1}u_j}\Bigr)\ee
and
\be\label{abi00f}{\cal T}_j(u_1, \dots ,u_j):=T_1(u_{j-1},u_j)\dots T_{j-1}(u_1,u_j)\ .\ee
Evaluating the numerical factors (\ref{abi00c}),
one obtains expressions for the symmetrizer and antisymmetrizer for the
algebra $BMW_n$ for $n>1$. The symmetrizer (its content sequence is $\{ 1,q^2,\dots ,q^{2(n-1)}\}$)
is equal to
\be\label{abi00ss}S_n=\frac{q^{n(n-1)/2}}{n_q!}\prod_{k=1}^{n-1}\frac{q^{2k-1}\nu^{-1}+1}{q^{4k-1}\nu^{-1}+1}
\ Y_n(1,q^2,\dots ,q^{2(n-1)})\ee
and the antiymmetrizer (its content sequence is $\{ 1,q^{-2},\dots ,q^{2(n-1)}\}$) is equal to:
\be\label{abi00aa}A_n=\frac{q^{-n(n-1)/2}}{n_q!}\prod_{k=1}^{n-1}\frac{q^{-2k-1}\nu^{-1}+1}{q^{-4k-1}\nu^{-1}+1}
\ Y_n(1,q^{-2},\dots ,q^{-2(n-1)})\ .\ee
Here we used the $q$-numbers, $n_q:=(q^n-q^{-n})/(q-q^{-1})$ and $n_q!:=2_q3_q\dots n_q$.

\vskip .2cm
The expressions (\ref{abi00ss}) and (\ref{abi00aa}) can be brought to a simpler form. Indeed, by the braid relation,
the product  $a_n:={\cal T}_n(1,\dots,q^{-2(n-1)})\dots {\cal T}_2(1,q^{-2})$
(this is the part of $Y_n(1,q^{-2},\dots ,q^{-2(n-1)})$ in (\ref{abi00aa})) can be rewritten as a word which
starts/ends by
$T_j(q^2)$ for any given $j=1,\dots,n-1$. Therefore $a_n T_i = T_i a_n = - q^{-1} a_n$ so
$a_n$ is proportional to $A_n$. Calculating $a_n^2$ (or, equivalently, evaluating the product of the operators
${\cal Q}$ in
(\ref{abi00aa}) on the antisymmetrizer) one finds the
proportionality coefficient and obtains the short factorized expression for the antisymmetrizer (\ref{abi00aa})
for the algebra $BMW_n$
which has the same form -- in terms of the
baxterized elements -- as the factorized expression  for the antisymmetrizer for the Hecke algebra
\be\label{abi00ab}A_n=\frac{(-1)^{n(n-1)/2}}{n_q!}\ {\cal T}_n(1,q^{-2},\dots ,q^{-2(n-1)})\dots {\cal T}_2(1,q^{-2})\ ,\ee
or
\be\label{abi00ab2}A_n=\frac{(-1)^{n-1}}{n_q}\ T_1(q^2)T_2(q^4)\dots T_{n-1}(q^{2(n-1)})\, A_{n-1}\ .\ee

The short factorized expression for the symmetrizer is obtained by applying the isomorphism
$BMW_n(-q^{-1},\nu)\to BMW_n(q,\nu)$ defined on generators by
$$BMW_n(-q^{-1},\nu)\ni T_j\mapsto T_j\in BMW_n(q,\nu)\ \ ,\ \ j=1,\dots,n-1\ .$$
The replacement $q\mapsto -q^{-1}$ in (\ref{a00}) and (\ref{a0}) produces the second set $T^*_i(u,v)$ and
$Q^*_i(u,v;c)$ of baxterized elements and converts (\ref{abi00ab}), (\ref{abi00ab2}) to
the short factorized expression for the symmetrizer (\ref{abi00ss})
\be\lb{abi00ab5} S_n = \displaystyle{
\frac{1}{n_q}\ T_1^*(q^{-2})T_2^*(q^{-4})\dots T_{n-1}^*(q^{-2(n-1)})\, S_{n-1}} \ee
 $$=\! \displaystyle{ \frac{1}{n_q!}\left[ T_1^*(q^{-2})T_2^*(q^{-4})\dots T_{n-1}^*(q^{-2(n-1)})\right]\!
\left[ T_1^*(q^{-2})T_2^*(q^{-4})\dots T_{n-2}^*(q^{-2(n-2)})\right] \cdots  T_1^*(q^{-2}) }.$$
The construction of all idempotents can be equivalently done with the help of the elements $T^*_i(u,v)$ and
$Q^*_i(u,v;c)$; the substitution $q\mapsto -q^{-1}$ into the expression for the $E_{U_n}$ leads to the idempotent
$E_{U_n^t}$ for the transposed up-down tableau.

The formulas (\ref{abi00ab2}) and (\ref{abi00ab5}) for the symmetrizer and antisymmetrizer were suggested in \cite{Isa}, see
also \cite{IsOg-BSS}. We refer the reader to the works \cite{Dipp}, \cite{Heck} and \cite{TuWe}, where
the symmetrizer and antisymmetrizer are given in several other forms.

\paragraph{5. Hecke algebra.} The fusion functions (\ref{abi00}) for the BMW algebra depend on $c$, whose value is
related to the parameter $\nu$ by (\ref{ccc}). In the Hecke algebra quotient of the BMW algebra (the quotient by the
ideal
generated by all $\kappa_i$) the parameter $\nu$ disappears. However the images of the fusion functions $\Psi_{U_n}$
in the
Hecke algebra quotient still contain non-trivially the parameter $c$. Considering $c$
as a free parameter in the Hecke algebra quotient we obtain a one parametric family of fusion functions (for the
Hecke
algebras), given by (\ref{abi00}), where we have to substitute
$$T_i(u_1,u_2)=T_i + \frac{q-q^{-1}}{u_2 / u_1-1}\; ,\;\;\; Q_i(u_1,u_2;c)=T_i+\frac{q-q^{-1}}{c u_1  u_2 -1}$$
in the definition (\ref{a4}). Explicitly the one-parametric family of the fusion functions for the Hecke algebra
reads
\be\lb{abi00h}
 \widetilde{\Psi}_{U_n}(u_1,u_2,\dots ,u_n)\! := \!
\left( \prod\limits_{k=1}^n \frac{u_k-c_k}{cu_kc_k-1}\frac{cu_k-1}{u_k-1}\right)
\widetilde{Y}_1(u_1)\widetilde{Y}_2(u_1,u_2)\dots\widetilde{Y}_n(u_1, \dots ,u_n) \, ,\ee
where
\be\lb{heck}\widetilde{Y}_j(u_1,\dots,u_j):= \prod^{\longleftarrow}_{k\colon j > k\geqslant 1}
\left(\! T_k + \frac{q-q^{-1}}{cu_ku_j-1}
\right)\prod^{\longrightarrow}_{k\colon 1\leqslant k < j}\left(\! T_k+\frac{q-q^{-1}}{u_j/u_k-1}\! \right)^{-1} ,\ee
$j=1,2,\dots$; the empty product is equal to 1, the symbol $\longleftarrow$ (respectively, $\longrightarrow$) over
$\prod$
means that the product is ordered according to the descend (respectively, ascend) of the product index $k$.

\vskip .2cm
The fusion function (\ref{abi00h}) generalizes the fusion function proposed in \cite{IsMoOs}; the results of
\cite{IsMoOs}
are reproduced at $c=0$ since $Q_i(u_1,u_2;0)=T_i^{-1}$.

\vskip .2cm
In the classical limit $q\rightarrow 1$, the family (\ref{abi00h}) of fusion functions leads to the one-parametric
family,
discussed in \cite{IMO}, of fusion functions for the symmetric group.

\paragraph{6. Reflection equation and fusion functions.} It was discovered in \cite{IsOg2} that the element
\be\lb{rea15}L_j(u)=\frac{cu\, y_{j} - 1}{u-y_{j}}\ee
solves the reflection equation
\be\lb{rea14}L_j(u)\, T_j\left(\frac{1}{cu\, v}\right)\,L_j(v) \,
T_j\left(\frac{u}{v}\right)=T_j\left(\frac{u}{v}\right)\, L_j(v)\,
T_j\left(\frac{1}{cu\, v}\right)\, L_j(u)\ee
in the BMW algebra. The element (\ref{rea15}) is just the same as in (\ref{a5}). This coincidence
is important in the study of the evaluation homomorphisms for the quantum universal enveloping algebras,
see \cite{IMO} for the classical counterpart.

\vskip .2cm
The main ingredients of the fusion function -- the elements ${\cal Y}_j$, $j=1,\dots,n-1$, defined
in (\ref{a4}) -- also satisfy the reflection equation
\be\lb{rea14b}\begin{array}{l}
{\cal Y}_j(u_1,\dots ,u_{j-1},v)\,
T_j\Bigl(\displaystyle{\frac{1}{cu\, v}}\Bigr)\, {\cal Y}_j(u_1,\dots ,u_{j-1},u) \,
T_j\Bigl(\displaystyle{\frac{u}{v}}\Bigr)^{-1}
 \\[1em]
\hspace{1cm}=  T_j\Bigl(\displaystyle{\frac{u}{v}}\Bigr)^{-1} \,{\cal Y}_j(u_1,\dots ,u_{j-1},u)\,
T_j\Bigl(\displaystyle{\frac{1}{cu\, v}}\Bigr)
\,{\cal Y}_j(u_1,\dots ,u_{j-1},v) \ .\end{array}\ee
This is shown by induction on $j$.

\vskip .2cm

The reflection equation (\ref{rea14b}) can be represented  graphically as

\unitlength=1.1cm
\begin{picture}(14,4)
\put(8.5,1){\line(1,0){5}}\put(8.5,0.97){\line(1,0){5}}
\put(7.7,3.4){\vector(1,-1){2.4}}\put(10.1,1){\vector(1,1){2}}
\put(7.6,3.15){\vector(3,-2){3.2}}\put(10.8,1){\vector(3,2){3}}
\put(9.5,0.7){\tiny ${\cal Y}_j(\vec{u},u)$}\put(10.7,0.7){\tiny ${\cal Y}_j(\vec{u},v)$}
\put(9.7,1.75){\tiny $T_{_j}\!(\frac{1}{c u v})$}\put(8.5,2.6){\tiny $T_j(\frac{u}{v})^{-1}$}
\put(7.4,2.8){\tiny $\frac{1}{cv}$}\put(7.9,3.4){\tiny $\frac{1}{cu}$}
\put(12,3.1){\tiny $u$}\put(13.5,3.1){\tiny $v$}\put(7.1,1.3){\Large $=$}
\put(1,1){\line(1,0){5}}\put(1,0.97){\line(1,0){5}}
\put(2.1,3){\vector(1,-1){2}}\put(4.1,1){\vector(1,1){2.5}}
\put(0.3,3){\vector(3,-2){3}}\put(3.3,1){\vector(3,2){3.5}}
\put(0,3.2){\tiny $\frac{1}{cv}$}\put(1.5,3.2){\tiny $\frac{1}{cu}$}
\put(6.2,3.4){\tiny $u$}\put(6.7,3.1){\tiny $v$}
\put(2.8,0.7){\tiny ${\cal Y}_j(\vec{u},v)$}\put(4.2,0.7){\tiny ${\cal Y}_j(\vec{u},u)$}
\put(3.3,1.8){\tiny $T_{_j}\!(\frac{1}{c u v})$}\put(5.7,2.3){\tiny $T_j(\frac{u}{v})^{-1}$}
\put(6,0.2){\bf Fig. 1}
\end{picture}

\noindent
Here we used the notation $\vec{u} := \{u_1,\dots,u_{j-1} \}$.
This picture gives the standard physical interpretation of the reflection
equation (\ref{rea14b}) which is used in the theory of the factorized scattering of
 particles on the half-line, see \cite{Che},\cite{Zam}.
The reflection from the boundary is described by the operator ${\cal Y}_j$ while the collision of
two particles is described by the operator $T_j$ or $T_j^{-1}$.
 The left hand side of the equality in Fig.~1 represents the process when: first,
the particle with the spectral parameter $\frac{1}{c v}$ reflects from the boundary and
changes its spectral parameter to $v$; second, this particle collides with
the particle having the spectral parameter $\frac{1}{c u}$; third, the particle with the
spectral parameter $\frac{1}{c u}$
reflects from the boundary and acquires the spectral parameter $u$; finally
the two reflected particles collide. Algebraically
this sequence of events is given by the left hand side of (\ref{rea14b}).
The picture in the right hand side of the equality in Fig.~1 shows the different scenario of the same process
(with different initial positions of the particles) and
algebraically is expressed by the right hand side of (\ref{rea14b}).
By the factorizability of the scattering both scenarios lead to the same result which is expressed
by the equation (\ref{rea14b}).

\vskip .2cm

One can write the fusion function (\ref{abi00}) in the form
(cf. (\ref{abi00b})):
\be\lb{abi00bb}
\! \Psi_{U_n}(u_1,u_2,\dots ,u_n) \! =\! \left( \prod\limits_{k=1}^n\frac{u_k-c_k}{cu_kc_k-1} \right)
 \mathfrak{Q}_n(u_1, \dots ,u_n)
 \, {\cal T}'_n(u_1,\dots,u_n) \dots {\cal T}'_2(u_1,u_2) \, ,
 \ee
where
\be\label{abi00ff}{\cal T}'_j(u_1, \dots ,u_j):=T_1(u_j,u_{j-1})^{-1}\dots T_{j-1}(u_j,u_1)^{-1} \ , \ee
\be\label{abi00dd}
 \mathfrak{Q}_n(u_1, \dots ,u_n):=
{\cal Y}_1(u_1) {\cal Q}_2(u_1,u_2) {\cal Y}_1(u_2)\dots {\cal Q}_n(u_1, ... ,u_n) {\cal Y}_1(u_n)  \; , \ee
with ${\cal Q}_j(u_1, \dots ,u_j)$ given by (\ref{abi00e}).

\vskip .2cm

Graphically, the part
\be
 \lb{abi00rr}
 \mathfrak{Q}_n(u_1, \dots ,u_n) \, {\cal T}'_n(u_1,\dots,u_n) \dots {\cal T}'_2(u_1,u_2) \; ,
 \ee
of the fusion function (\ref{abi00bb}) is represented as the process of scattering
of $n$ particles on the half-line.
For example, the expression (\ref{abi00rr}) for $n=4$
 \be\label{abi00d4}
 \begin{array}{c}
{\cal Y}_1(u_1) {\cal Q}_2(u_1,u_2) {\cal Y}_1(u_2) {\cal Q}_3(u_1, u_2 ,u_3)\
 {\cal Y}_1(u_3) {\cal Q}_4(u_1, u_2 ,u_3 ,u_4)\  {\cal Y}_1(u_4) \cdot \\ [0.2cm]
  {\cal T}'_4(u_1,u_2,u_3,,u_4) \, {\cal T}'_3(u_1,u_2,u_3) \, {\cal T}_2'(u_1,u_2)  \; ,
\end{array}
 \ee
is visualized in Fig.~2:

\noindent
\unitlength=1.2cm\begin{picture}(12,6)
\put(2,1){\line(1,0){8}}\put(2,0.96){\line(1,0){8}}\put(6,3.1){\vector(1,-3){0.7}}
 \put(6.7,1){\vector(1,3){1.3}}\put(3.9,3){\vector(1,-1){2}}
\put(5.9,1){\vector(1,1){4}}\put(1.8,3){\vector(3,-2){3}}\put(4.8,1){\vector(3,2){6}}
 \put(0.4,2.4){\vector(2,-1){2.8}}\put(3.2,1){\vector(2,1){8}}
\put(0,2.5){\tiny $\frac{1}{cu_1}$}\put(1.7,3.2){\tiny $\frac{1}{cu_2}$}\put(3.5,3.2){\tiny
$\frac{1}{cu_3}$}\put(5.5,3.2){\tiny $\frac{1}{cu_4}$}
\put(11.2,4.8){\tiny $u_1$}\put(9.3,4.8){\tiny $u_3$}\put(7.5,4.5){\tiny $u_4$}\put(10.3,5){\tiny $u_2$}
\put(2.5,0.7){\tiny ${\cal Y}_1(u_1)$}\put(4.3,0.7){\tiny ${\cal Y}_1(u_2)$}\put(5.5,0.7){\tiny ${\cal
Y}_1(u_3)$}\put(6.7,0.7){\tiny ${\cal Y}_1(u_4)$}
\put(4,1.6){\tiny $Q_1$}\put(4.5,1.9){\tiny $Q_2$}\put(5.8,2.6){\tiny $Q_3$}\put(5,1.4){\tiny $Q_1$}
\put(6.2,1.6){\tiny $Q_1$}\put(6,2.1){\tiny $Q_2$}
\put(7.1,2){\tiny $T^{-1}_1$}\put(6.85,2.65){\tiny $T^{-1}_2$}\put(7,3.2){\tiny $T^{-1}_3$}
\put(8.2,3.1){\tiny $T^{-1}_1$}\put(8.3,3.8){\tiny $T^{-1}_2$}\put(9.4,3.9){\tiny $T^{-1}_1$}\put(6,0){\bf Fig. 2}
\end{picture}
\vspace{0.5cm}

\noindent
In Fig. 2 the spectral parameters for operators $T^{-1}_k$ and $Q_k$  can be restored by the rules:

\vspace{0.5cm}
\unitlength=8mm\begin{picture}(19,2)
\put(0.5,0.8){$T_{i-j}(u_i,u_j)^{-1} \;\; =$}
 \put(4.3,0.1){$j$}\put(6.5,0.1){$i$}\put(3.9,1.7){$u_i$}\put(6.7,1.7){$u_j$}\put(4.7,0.2){\vector(1,1){1.5}}
\put(6.2,0.2){\vector(-1,1){1.5}}\put(7.4,0.8){,}\put(8.1,0.8){$Q_{j-i}(u_i,u_j;c) = T_{j-i}(\frac{1}{cu_j},u_i) \;
=$}
\put(14.6,0.1){$j$}\put(16.8,0.1){$i$}\put(14.4,1.7){$u_i$}
 \put(16.7,1.7){$\frac{1}{cu_j}$}\put(15,0.2){\vector(1,1){1.5}}\put(16.5,0.2){\vector(-1,1){1.5}}
\put(17.4,0.8){.}
\end{picture}

\noindent
{\bf Remark.}
There is another interpretation of the expressions (\ref{abi00rr}) for different $n$.
Namely these are fusion solutions of the so-called higher reflection equations.
The fusion method
 was proposed in \cite{KRS},\cite{KS}
for the Yang--Baxter equation and in \cite{MN}
for the reflection equation. The interpretation of (\ref{abi00rr}) as a fusion solution of
reflection equation justifies the name ``fusion function" for $\Psi_{U_n}$ in (\ref{abi00}).

\section{Classical limit}\label{secclli}
\setcounter{equation}{0}

In this section we demonstrate that the fusion procedures of \cite{IsMol} and \cite{IMO}
for the Brauer algebra can be obtained from the fusion procedure for the
BMW algebra
by taking two classical limits. Following \cite{IMO}, we denote the generators
of the Brauer algebra ${\cal B}_n(\omega)$ by $s_1,\dots,s_{n-1},\epsilon_1,\dots,\epsilon_{n-1}$.  However, we denote
by
$\theta_i$ the rational
spectral parameters instead of the $u_i$ used in \cite{IMO} to avoid confusion.

\vskip .2cm
The symbol $\to$ in this Section means ``tends to".

\paragraph{1. Contraction of the BMW algebra to the Brauer algebra.} The BMW algebra is a flat deformation of the
Brauer algebra;
expressed differently, the Brauer algebra is a contraction of the BMW algebra. To describe the contraction, we set
\be\lb{sc1} q = e^h\; ,\;\;\; \nu =q^{1-\omega}=e^{h(1-\omega)}\; \ee
and let $h$ tend to 0,
\be\lb{sc0b}h\to 0\ .\ee
In the contraction, the generator $T_i$ tends to $s_i$ which has order 2; in the limit (\ref{sc0b}) the element
$\kappa_i$ becomes thus an independent generator; we have
\be\lb{sc0}T_i\to s_i\; ,\;\;\; \kappa_i \to \epsilon_i\; .\ee
The expression (\ref{sc1}) for $\nu$ contains $\omega$ which becomes the parameter of the Brauer algebra family.

\paragraph{2. Behavior of the spectral parameters.} The description of the contraction is given by the
formulas (\ref{sc1}), (\ref{sc0b}) and (\ref{sc0}). However
the fusion function for the BMW algebra depends on the spectral parameters $u_j$ and we have to complete the
formulas  (\ref{sc1}), (\ref{sc0b}) and (\ref{sc0}) for the contraction with the information about the behavior of the
spectral
parameters; the relation to the rational spectral parameters is given by
\be\lb{sc2}u_j=e^{2h(\theta_j-\frac{\omega-1}{2})}\ee
for  the classical limit, determined by (\ref{sc1}), (\ref{sc0b}) and (\ref{sc0}).

\vskip .2cm
To explain the prescription (\ref{sc2}) recall the
definition of the classical content of boxes in an up-down tableau $U_n=(\Lambda_1,\Lambda_2,\dots,\Lambda_n)$. {}For
the box in the $a$-th row and $b$-th column the classical content is defined by
\be\label{doqcfudte}(b-a)+\frac{\omega-1}{2}\ee
if this box was added to the diagram $\Lambda_{j-1}$ to obtain the diagram $\Lambda_j$ and
\be\label{doqcfudtf}-(b-a)-\frac{\omega-1}{2}\ \ee
if this box was removed from the diagram $\Lambda_{j-1}$ to obtain the diagram $\Lambda_j$.

\vskip .2cm
The quantum content is defined in (\ref{doqcfudt}) and (\ref{doqcfudtb}). In (\ref{sc2}) we made the shift in the
exponent
for the agreement between the quantum and classical contents.

\paragraph{3. Classical limit of the fusion function.} Performing the classical limit, determined by (\ref{sc0b}),
(\ref{sc1}),
(\ref{sc0}) and (\ref{sc2}), of the building
blocks $Q(u,v)$ and $T(u,v)$ of the fusion function, we find
\be\lb{sc3}\begin{array}{c}
\displaystyle{Q_i(u_1,u_2)\to s_i-\frac{\epsilon_{i}}{\theta_1+\theta_2}
=s_i\Bigl(1-\frac{\epsilon_{i}}{\theta_1+\theta_2}\Bigr)\ ,}\\[1.5em]
\displaystyle{T_i(u_1,u_2)\to s_{i}-\frac{1}{\theta_1-\theta_2}
=s_i\Bigl(1-\frac{s_i}{\theta_1-\theta_2}\Bigr)}\ ,\end{array}\ee
so this classical limit of (\ref{abi00}) reproduces the fusion
function from \cite{IsMol}.

\paragraph{4. Second contraction.} There exists a different contraction of the BMW algebra to the Brauer algebra. It
is determined by (\ref{sc0b}), (\ref{sc0}), and
\be\lb{sc1b} q =- e^h\; ,\;\;\; \nu =e^{h(\omega-1)}\; \ee
instead of (\ref{sc1}).

\paragraph{5. Behavior of the spectral parameters for the second contraction.} As we shall see, in the limit
determined
by  (\ref{sc0b}), (\ref{sc1b}) and (\ref{sc0}), we obtain the
idempotents for the transposed diagrams and up-down tableaux. This phenomenon is clearly illustrated on the
example of the (anti)symmetrizers. Indeed, on the $q$-deformed level, the symmetrizer is the projector on the
eigenspace
with the eigenvalue $q$. With the prescription (\ref{sc0b}) the parameter $q$ tends to (-1) so the $q$-symmetrizer
becomes the classical antisymmetrizer in this limit.

\vskip .2cm
With the prescription (\ref{sc1b}), the quantum and classical contents cannot be made compatible within the Ansatz
$u=e^{hf(\theta)}$ with any function $f$ ($u$ is the quantum content, $\theta$ is the classical one); the definition
of the
classical content of a box has to be modified. {}Following the above example of the (anti)symmetrizers, we define the
t-classical content of a box to be the usual classical content of the transposed box (``t" in t-classical stands for
transposition).
We shall write explicitly the analogues of the formulas (\ref{doqcfudte}) and (\ref{doqcfudtf}) for an up-down
tableau
$U_n=(\Lambda_1,\Lambda_2,\dots,\Lambda_n)$. {}For the box in the $a$-th row and $b$-th column the t-classical
content is defined by
\be\label{doqcfudtc}-(b-a)+\frac{\omega-1}{2}\ee
if this box was added to the diagram $\Lambda_{j-1}$ to obtain the diagram $\Lambda_j$ and
\be\label{doqcfudtd}(b-a)-\frac{\omega-1}{2}\ \ee
if this box was removed from the diagram $\Lambda_{j-1}$ to obtain the diagram $\Lambda_j$.

\vskip .2cm
Instead of (\ref{sc2}) we impose now
\be\lb{sc2b}u_j=e^{2h(-\theta_j+\frac{\omega-1}{2})}\ee
for the agreement of the quantum and t-classical contents.

\paragraph{6. Second classical limit of the fusion function.} Performing the classical limit, determined by
(\ref{sc0b}),
(\ref{sc1b}), (\ref{sc0}) and (\ref{sc2b}), of the elements $Q(u,v)$ and $T(u,v)$, we find
\be\lb{sc3b}\begin{array}{c}
\displaystyle{Q_i(u_1,u_2)\to s_i+\frac{1}{\theta_1+\theta_2-\varkappa}-\frac{\epsilon_{i}}{\theta_1+\theta_2}
=s_i\Bigl(1+\frac{s_i}{\theta_1+\theta_2-\varkappa}-\frac{\epsilon_{i}}{\theta_1+\theta_2}\Bigr)\ ,}\\[1.5em]
\displaystyle{T_i(u_1,u_2)\to s_{i}-\frac{1}{\theta_1-\theta_2}+\frac{e_i}{\theta_1-\theta_2-\varkappa}
=s_i\Bigl(1-\frac{s_i}{\theta_1-\theta_2}+\frac{e_i}{\theta_1-\theta_2-\varkappa}\Bigr)}\ ,\end{array}\ee
where $\varkappa=\frac{\omega}{2}-1$. Thus the function (\ref{abi00}) in this limit tends precisely to the fusion
function introduced in \cite{IMO}.

\vskip .2cm
We repeat that the evaluation of the fusion function (\ref{abi00}) for the BMW algebra on contents of an up-down
tableau $U_n$ descends to the evaluation of the limiting fusion function for the Brauer algebra on the transposed
tableau $U_n^T$. The idempotents belonging to the diagram $\lambda$ of the BMW algebra tend to the idempotents
belonging to the diagram $\lambda^T$ of the Brauer algebra.

\paragraph{7. Further contractions.} In the two contractions above, $q$ was tending to 1 or (-1). When we fix the
behavior of $q$ there are
two possibilities for the behavior of $\nu$: (\ref{lammu}), considered as an equation for $\nu$, is quadratic and we
must have $\mu\rightarrow\omega$ in the classical limit. Therefore, the BMW-algebra contracts to the Brauer algebra
in two more regimes ($h\to 0 $ as before):
\be\lb{sc4b}q=e^h\; ,\;\;\;\nu =-e^{h(\omega-1)}\ee
and
\be\lb{sc1c} q =- e^h\; ,\;\;\; \nu =-e^{h(1-\omega)}\ .\ee
{}For these regimes we have to impose that the generator $T_i$ tends to $(-s_i)$ in the classical limit to respect the
relations
$s_i\epsilon_i=\epsilon_is_i=\epsilon_i$ of the Brauer algebra. Thus, (\ref{sc0}) gets replaced by
\be\lb{sc5}T_i\to -s_i\; ,\;\;\;\kappa_i\to \epsilon_i\; .\ee
We omit further details of the corresponding limits of the fusion function (3.16).


\begin{thebibliography}{99}
\bibitem{BeBla} A. Beliakova and C. Blanchet, {\it Skein construction of idempotents in Birman--Murakami--Wenzl
algebras}.  Math. Ann.  {\bf 321}  no. 2  (2001) 347 -- 373.
\bibitem{BirWen}  J. Birman and H. Wenzl, {\it Braids, link polynomials and
a new algebra}. Trans. Amer. Math. Soc. {\bf 313} (1989) 249--273.
\bibitem{CGX} Y. Cheng, M. L. Ge and K. Xue, {\it Yang--Baxterization of braid group representations}.
Comm. Math. Phys.  {\bf 136} no. 1 (1991) 195--208.
\bibitem{Cher1} I. V. Cherednik, {\it A new interpretation of Gelfand--Tzetlin bases}.
Duke Math. J. {\bf 54} no. 2 (1987) 563--577.
\bibitem{Che} I. V. Cherednik, {\it Factorizing particles on a half-line and root systems}.
Theor. Math. Phys., {\bf 61} no. 1 (1984) 977--983.
\bibitem{Dipp} R. Dipper, Jun Hu, F. Stoll, {\it Symmetrizers and antisymmetrizers for
the BMW algebra}. arXiv: 1109.0342.
\bibitem{Heck} I. Heckenberger, A. Sch\"{u}ler, {\it Symmetrizer and antisymmetrizer of the
Birman-Wenzl-Murakami algebras}. Lett. Math. Phys., {\bf 50} (1) (1999) 45--51.
\bibitem{Zam} S. Ghoshal, A. B. Zamolodchikov, {\it Boundary S-Matrix and boundary state in two-dimensional integrable
    quantum field theory}.
Int. J. Mod. Phys. A{\bf 9} (1994) 3841--3886; arXiv: hep-th/9306002
\bibitem{GHJ} F. M. Goodman, P. de la Harpe and V. F. R. Jones, {\it Coxeter graphs and towers of
algebras}, Springer (1989).
\bibitem{Isa} A. P. Isaev, {\it Quantum groups and Yang-Baxter equations}.
preprint MPIM (Bonn), MPI 2004-132 (2004), \\
 http://webdoc.sub.gwdg.de/ebook/serien/e/mpi$\underline{\phantom{a}}$mathematik/2004/132.pdf
\bibitem{IsMol} A. P. Isaev and A. I. Molev, {\it Fusion procedure for the Brauer algebra}.
Algebra i Analiz {\bf 22} no. 3 (2010) 142--154; arXiv: 0812.4113 [math.RT]
\bibitem{IMO} A. P. Isaev, A. I. Molev and O. V. Ogievetsky, {\it A new fusion procedure for the Brauer
algebra and evaluation homomorphisms}. Int. Math. Res. Not. 2011 doi: 10.1093/imrn/rnr126; arXiv: 1101.1336
[math.RT]
\bibitem{IsMoOs} A. P. Isaev, A. I. Molev and A. F. Os'kin, {\it On the idempotents of Hecke algebra}.
Lett. Math. Phys. {\bf 85} (2008) 79-90; arXiv: 0804.421 [math.QA]
\bibitem{IsOg3} A. P. Isaev and O. V. Ogievetsky, {\it On representations of Hecke algebras}. Czech.
Journ. Phys. {\bf 55} no. 11 (2005) 1433--1441.
\bibitem{IsOg2} A. P. Isaev and O. V. Ogievetsky, {\it On Baxterized solutions of reflection equation and
integrable chain models}. Nucl. Phys. B {\bf 760} (2007) 167--183; arXiv: math-ph/0510078
\bibitem{IsOg-BSS} A. P. Isaev and O. V. Ogievetsky, {\it Braids, shuffles and symmetrizers.} J. Phys. A:
Math. Theor. {\bf 42} (2009) 1-15; arXiv: math.QA/0511618
\bibitem{IsOg} A. P. Isaev and O. Ogievetsky, {\it Jucys--Murphy elements for Birman--Murakami--Wenzl
algebras}. Physics of Particles and Nuclei Letters {\bf 8} no. 3  (2011) 394--407; arXiv: 0912.4010 [math.QA]
\bibitem{Jones} V. F. R. Jones, {\it On a certain value of the Kauffman polynomial}.
Comm. Math. Phys. {\bf 125} (1989) 459--467.
\bibitem{Ju0} A. Jucys, {\it On the Young operators of the symmetric group}. Lietuvos Fizikos Rinkinys {\bf 6} (1966)
    163--180.
\bibitem{Ju} A. Jucys, {\it Factorization of Young projection operators for the symmetric group}. Lietuvos
Fizikos Rinkinys {\bf 11} (1971) 5--10.
\bibitem{KRS} P. P. Kulish, N. Yu. Reshetikhin and E. K. Sklyanin, {\it Yang--Baxter equation and representation
    theory: I}.
Lett. in Math. Phys. {\bf 5} no. 5 (1981) 393--403.
\bibitem{KS} P. P. Kulish and E. K. Sklyanin, {\it Quantum spectral transform method: recent developments.} In:
    Integrable Quantum Field Theories,
Lect. Notes Phys. {\bf 151} (1982) 61--119.
\bibitem{LeRa} R. Leduc and A. Ram, {\it A ribbon Hopf algebra approach to the irreducible representations of
    centralizer algebras: The Brauer,
Birman--Wenzl, and type A Iwahori-Hecke algebras}.
Advances in Mathematics {\bf 125} no. 1 (1997) 1-- 94.
\bibitem{MN} L. Mezinchescu and R. Nepomechie, {\it Fusion procedure for open chains}. J. Phys. A: Math. Gen. {\bf 25}
    (1992) 2533--2543.
\bibitem{Mol} A. I. Molev, {\it On the fusion procedure for the symmetric group}. Rep. Math. Phys. {\bf 61} (2008)
181--188; arXiv: math/0612207 [math.RT]
\bibitem{Mur00} J. Murakami,  {\it The Kauffman polynomial of links and representation theory}.
Osaka J. Math. {\bf 24} (1987) 745--758.
\bibitem{Mur0} J. Murakami, {\it The representation of the q-analogue of Brauer's centralizer algebras
and the Kauffman polynomial of links}. Publ. RIMS {\bf 26} no. 6 (1990) 935--945.
\bibitem{Mur} J. Murakami, {\it Solvable lattice models and algebras of face operators}. Adv. Studies in Pure
Math. {\bf 19} (1989) 399--415.
\bibitem{Mu} G. E. Murphy, {\it A new construction of Young's seminormal representation of the symmetric
groups}. J. Algebra {\bf 69} (1981) 287--297.
\bibitem{Naz} M. Nazarov, {\it Yangians and Capelli identities}. In: Kirillov's Seminar on Representation Theory (G.
    I.  Olshanski, ed.), {Amer. Math. Soc. Transl.} {\bf 181}, Amer. Math. Soc.,
Providence, RI  (1998) 139--163.
\bibitem{Naz2} M. Nazarov, {\it A mixed hook-length formula for affine Hecke algebras}.
European J. Combin. {\bf 25} (2004) 1345--1376.
\bibitem{O} O. V. Ogievetsky, {\it Uses of quantum spaces}.
Contemp. Math. {\bf 294}, Amer. Math. Soc., Providence, RI (2002) 161-232.
\bibitem{OPdA} O. Ogievetsky and L. Poulain d'Andecy, {\it On representations of cyclotomic Hecke algebras}.
Mod. Phys. Lett. A {\bf 26} no. 11 (2011) 795--803; arXiv: 1012.5844 [math-ph]
\bibitem{OV}A. Okounkov and A. Vershik, {\it A new approach to representation theory of symmetric groups
II}. Selecta Math (New series) {\bf 2} no. 4 (1996) 581--605.
\bibitem{RaWe} A. Ram and H. Wenzl, {\it Matrix units for centralizer algebras}. J. Algebra {\bf 145} no. 2 (1992)
 378--395.
 \bibitem{TuWe} I. Tuba, H. Wenzl, {\it On braided tensor categories of type $BCD$}.
 J. Reine Angew. Math.  {\bf 581}  (2005) 31--69.
 \bibitem{We} H. Wenzl, {\it Quantum groups and subfactors of type B, C and D}. Comm. Math. Phys. {\bf 133} (1990)
    383--432.

\end{thebibliography}
\end{document}